\documentclass[12pt]{article}
\usepackage{amssymb}
\usepackage{epsf}
\usepackage{amsmath}
\usepackage{mathptm}
\usepackage{pst-plot}
\usepackage{pst-grad}
\usepackage{pst-node}
\usepackage{amsmath}
\input amssym.def
\input amssym.tex

\def \Z{{\Bbb Z}}

\newtheorem{Th}{THEOREM}[section]
\newtheorem{dfn}[Th]{DEFINITION}
\newtheorem{prop}[Th]{PROPOSITION}

\newtheorem{lemma}[Th]{LEMMA}

\begin{document}

\title{$K$-TWISTED EQUIVARIANT $K$-THEORY FOR $SU(N)$}
\author{Bin Zhang \\
Department of Mathematics\\
The State University of New York\\
Stony Brook, NY 11794-3651\\
\textit{bzhang@math.sunysb.edu} \thanks{%
The author would like to thank the Max-Planck-Institut f\"ur Mathematik at
Bonn for providing me the wonderful working environment and financial
support during my visit at MPIM 06/06/2003-09/01/2003, where this work is
finally done.} }
\date{}
\maketitle

\begin{abstract}
We present a version of twisted equivariant $K$-theory-$K$-twisted
equivariant $K$-theory, and use Grothendieck differentials to compute the $K$
-twisted equivariant $K$-theory of simple simply connected Lie groups. We
did the calculation explicitly for $SU(N)$ explicitly. The basic idea is to interpret an equivariant gerbe as an element of equivariant $K$-theory of degree 1.
\end{abstract}

\section{Introduction}

Let $G$ be a finite dimensional simple Lie group, a classical question
related to it is to understand the space $Hom(\pi ,G)/G$, where $\pi $ is a
finitely presented group. This space $Hom(\pi ,G)/G$ is the moduli space of
flat connections on a principal $G$-bundle on a manifold with fundamental
group $\pi $. Because of Atiyah and Segal's result \cite{AS1}, and the fact
that $K$-theory is defined for a large class of geometric objects including
usual topological spaces and non-commutative ones, our first approach is to
study the equivariant $K$-theory of $Hom(\pi ,G)$. We get the answer for the
case $\pi =\mathbb{Z}$, i.e. the equivariant $K$-theory $K_{G}^{\ast
}(G)\cong \Omega _{R(G)/\mathbb{Z}}^{\ast }$ \cite{BZ1} the algebra of
Grothendieck differentials of the representation ring $R(G)$ of $G$ over $
\mathbb{Z}$ when $G$ is compact and the fundamental group is torsion free
(the general situation is still open). This is the origin of our project
about Grothendieck differentials in $K$-theory.

We get interested in twisted $K$-theory because of Freed-Hopkins-Teleman's
result on twisted equivariant $K$-theory and Verlinde algebras \cite{Fd1}, 
\cite{Fd2}, \cite{FHT1}, that is for a Lie group $G$, the Verlinde algebra $
V_{k}(G)$ at level $k$ is twisted equivariant $K$-theory of $G$ (with
adjoint action) at particular degree. Unfortunately, they didn't publish
their proof yet. The main idea for this paper is to use Grothendieck
differentials to give a partial proof of their result, and supply a candidate for the geometric definition of twisted $K$-theory.

The first question we need solve is to find a good geometric model for
twisted equivariant $K$-theory. Let $\mathbb{H}$ be a infinite dimensional
separable Hilbert space, and $U=U(\mathbb{H})$ be the set of unitary
operators on $\mathbb{H}$, we know that $U$ is contractible. The group $U$
has a natural subgroup $\{e^{i\theta }I\}$ which is isomorphic to $S^{1}$,
let us denote the quotient group by $PU$. For a topological space $X$, in
principle, a twistor is a principal $PU$-bundle over $X$, thus an element in 
$H^{3}(X,\mathbb{Z})$. Naturally a geometric realization of $H^{3}(X,\mathbb{
Z})$ elements is needed. We already have a geometric realization of $H^{3}$
classes, i.e., gerbes \cite{Bj1}. Based on the idea of gerbes, there are
some other geometric realizations, like bundle gerbes \cite{Mm1} or central
extensions of groupoids \cite{BXZ1} \cite{TXL1}. All these involves infinite
dimensional objects. We are more interested in finite dimensional
realization of gerbes, like \cite{Me1} \cite{GR1}. But how can we do
twisting with gerbes? As far as I know, there is no clean geometric
definition for twisted $K$-theory. The equivariant situation is more subtle,
in this case, whether to use equivariant gerbes \cite{Bj2} to do twisting is
questionable.

We present a solution to these questions in nice situation. We study the
twisted equivariant $K$-theory of $G$ (with adjoint action). In this case,
we interpret an equivariant gerbe as an element of $K^1$, then based on this 
$K^1$ element, we give an intuitive definition of $K$-twisted equivariant $K$
-theory. The paper is basically two parts. In the first part, we prove that
an element in the equivariant cohomology $H^3_G(G)$ can be interpreted as
an element of $K^1_G(G)$, and in the second part, we use the definition we
give to do calculation for $SU(N)$ explicitly (in fact we did the
calculation for classical groups, but for simplity and to demonstrate the
idea, we just present the case for $SU(N)$).

\section{$K$-Twisted $K$-theory}

In the section, we first discuss the general picture of twisted $K$-theory
and then present our definition for $K$-twisted (equivariant) $K$-theory.

$K$-theory is a generalized cohomology theory \cite{Am1}. For a paracompact
topological space $X$, $K^*(X)$ has several equivalent definitions:

1. Geometric definition: equivalence classes of complex of vector bundles
over $X$.

2. Homotopic definiton: Homotopy classes of maps: $[X,Fred]$, $[X,Fred_{as}]$
, where $Fred$ and $Fred_{as}$ are the set of Fredholm operators and
self-adjoint operators in $\mathbb{H}$.

3. Algebraic definition: $K$-theory of $C^{\ast }$-algebra $C_{0}(X)$. 

Based on the homotopic definition of $K$-theory, the general picture of the
twisted $K$-theory can be as follows. If we have a principal $PU$-bundle $P$
over $X$, notice there are natural actions of $PU$ on $Fred$ and $Fred_{as}$
, we can form the spaces $P\times _{PU}Fred=(P\times Fred)/PU$ and $P\times
_{PU}Fred_{as}$, which are fiber bundles over $X$, then we can define the
twisted $K$-theory as the homotopy classes of sections of these two bundles.
There are general definitions of twisted $K$-theory from point of view of $
C^{\ast }$-algebra, see \cite{Rj1}, or \cite{TXL1} for the equivariant cases
for detail. We are more interested in a geometric picture of twisted $K$
-theory, and if possible, a definition with finite dimensional objects.

The twistor, i.e., the principal $PU$-bundle over $X$ is classified by $
H^{1}(X,PU)$. The exact sequence of groups $1\rightarrow S^{1}\rightarrow
U\rightarrow PU\rightarrow 1$ implies that $PU$ is a model for $BS^{1}$, the
classifying space of $S^{1}$. Thus $H^{1}(X,PU)\cong H^{2}(X,S^{1})\cong
H^{3}(X,\mathbb{Z})$, So the twistor is classified by $H^{3}(X,\mathbb{Z})$.
The geometric construction of a class in $H^{3}(X,\mathbb{Z})$ is a gerbe 
\cite{Bj1}. In brief, we use a gerbe to do twisted $K$-theory.

One might hope to use vector bundles to construct the twisted $K$-theory
geometrically. This is succeeded only in case that the twistor is a torsion
element in $H^{3}(X,\mathbb{Z})$ \cite{BCMKS}. In this case the twisted $K$
-theory is the Grothendieck group of the category of twisted bundles. The
essential problem is the non-existence of finite dimensional twisted bundles
in general.

The geometric picture for the twisted equivariant $K$-theory is more subtle.
Let $G$ be a topological group, $X$ be a $G$-space, the equivariant $K$
-theory $K_{G}^{\ast }(X)$ can be defined in the similar ways \cite{Sg1}.
The question in this case is what kind of twistor we can use. The natural
generalization of non-equivariant case is the elements in $H_{G}^{3}(X)$,
the 3rd degree equivariant cohomology, in other words equivariant gerbes.
But there is some problem if we use it to a geometric approach. The reason
is that an element of $H_{G}^{3}(X)$ is an object on $EG\times _{G}X$, not
exactly an equivariant object on $X$. There is a question just like
non-equivariant case,  what kind geometric objects we can use, again the
non-existence of twisted equivariant bundle is a problem.

There is another point of view for the whole picture. Let $X$ be a finite
dimensional object, for example, a finite dimensional manifold, then the
chern character $ch:K^{1}(X)\otimes \mathbb{Q}\cong H^{odd}(X,\mathbb{Q})$.
So up to $\mathbb{Z}$-torsion, an element in $H^{3}(X,\mathbb{Z})$ can be
viewed as an element in $K^{1}(X)$. This simple observation suggests the
following intuitive definition of $K$-twisted $K$-theory.

\begin{dfn}
Let $X$ be a topological space, and $\alpha \in K^{1}(X)$, the $K$-twisted $K
$-theory $^{\alpha }K^{\ast }(X)$ is the homology of the following complex,

\begin{equation*}
\cdots \overset{\wedge \alpha}{\to} K^0(X)\overset{\wedge \alpha}{\to} K^1(X)
\overset{\wedge \alpha}{\to} K^0(X)\overset{\wedge \alpha}{\to}\cdots
\end{equation*}
\end{dfn}

The desired properties of twisted $K$-theory are obvious from this
Definition. This definition should agree with the homotopic definition in case 
$\alpha $ is a non-torsion element in $H^{3}(X,\mathbb{Z})$, and there
should be a more general geometric definition of twisted $K$-theory which
generalizes this definition and twisted bundle in the torsion case. We are
working on this topic.

This definition can be easily generalized to the equivariant case,

\begin{dfn}
Let $X$ be a topological space, $G$ be a compact topological group acting on 
$X$, and $\alpha \in K_{G}^{1}(X)$, the $K$-twisted $K$-theory $^{\alpha
}K_{G}^{\ast }(X)$ is the homology of the following complex,

\begin{equation*}
\cdots \overset{\wedge \alpha}{\to} K^0_G(X)\overset{\wedge \alpha}{\to}
K^1_G(X)\overset{\wedge \alpha}{\to} K^0_G(X)\overset{\wedge \alpha}{\to}
\cdots
\end{equation*}
\end{dfn}

\section{The basic gerbe as an element of $K^1_G(G)$}

Let $G$ be a $n$-dimensional compact simple simply-connected Lie group of
rank $d$, $T$ be a maximal torus of $G$, and $W$ be the Weyl group of $G$
with respect to $T$. We use $R(G)$, $R(T)$ to denote the representation
rings of $G$ and $T$ respectively. If $\chi _{1},\chi _{2},\cdots ,\chi _{d}$
are the simple characters of $T$, then the character group $X^{\ast
}(T)=Hom(T,S^{1})$ is the free abelian group generated by $\chi _{1},\chi
_{2},\cdots ,\chi _{d}$, and the representation ring $R(T)$ is the group
ring $\mathbb{Z}[X^{\ast }(T)]=\mathbb{Z}[\chi _{1},\chi _{2},\cdots ,\chi
_{d},\chi _{1}^{-1},\chi _{2}^{-1},\cdots ,\chi _{d}^{-1}]$. The Weyl group $W$ acts on $R(T)$, the invariant subalgebra $R(T)^{W}$ is the
representation ring $R(G)$, which is a polynomial ring generated by
\textquotedblleft basic\textquotedblright\ representations $\rho _{1},\rho
_{2},\cdots ,\rho _{d}$ corresponding to a choice of a set of simple roots.

The cohomology of $T$ can be easily described in terms of these characters.
The character $\chi _{i}:T\rightarrow S^{1}$ can be viewed as an element of $
[X,S^{1}]\cong H^{0}(X,S^{1})\cong H^{1}(X,\mathbb{Z})$, let us denote this
element by $\eta _{i}$. By this way, we get a homomorphism of abelian groups 
$X^{\ast }(T)\rightarrow H^{1}(T,\mathbb{Z})$, and $H^{\ast }(T,\mathbb{Z}
)\cong \bigwedge (\eta _{1},\eta _{2},\cdots ,\eta _{d})$. 

The $K$-theory can be described in similar way. A character $\chi
_{i}:T\rightarrow S^{1}=U(1)$ defines a line bundle over the suspension of $T
$, thus defines an element of $K^{1}(T)$, again we denote this element by $
\eta _{i}$. Therefore we have a homomorphism between abelian groups: $
X^{\ast }(T)\rightarrow K^{1}(T)$, and $K^{\ast }(T)\cong \bigwedge (\eta
_{1},\eta _{2},\cdots ,\eta _{d})$. In particular we see that there is an
isomorphism $c:K^{\ast }(T)\cong H^{\ast }(T,\mathbb{Z})$, where the map is
in fact the first chern class of bundles, and this map is equivariant under
the action of Weyl group $W$.

Let $X$ be a paracompact space, $H$ be a compact topological group acting on 
$X$, then the equivariant cohomology is defined as $H_{H}^{\ast }(X)=H^{\ast
}(EH\times _{H}X)$ \cite{Ba1}, where $EH\rightarrow BH$ is a universal
principal $H$-bundle, $BH$ is a classifying space for $H$. In particular, $%
H_{H}^{\ast }(pt)=H^{\ast }(BH)$, and the bundle map $EH\times
_{H}X\rightarrow BH$ give $H_{H}^{\ast }(X)$ a $H_{H}^{\ast }(pt)$-module
structure.

In the case of the torus $T$, the coefficient ring $H^{\ast }(BT)$ can also
be described in terms of the character group $X^{\ast }(T)$. For any
character $\chi :T\rightarrow S^{1}$, it defines a line bundle $ET\times _{T}%
\mathbb{C}\chi $ over $BT$, the first chern class of this bundle gives an
abelian group homomorphism: $X^{\ast }(T)\rightarrow H^{2}(BT)$, this
induces an isomorphism between $H^{\ast }(BT)$ and the symmetric algebra $%
S_{T}$ of $X^{\ast }(T)$. Notice that $H^{\ast }(BT)$ carries a natural
action of the Weyl group $W$.

Let us consider $G$ as a $G$-space with adjoint action, it is well-known
that $H_{G}^{3}(G,\mathbb{Z})\cong \mathbb{Z}$, and the generator (up to
sign) is called the basic (equivariant) gerbe. There are several ways to
describe this gerbe \cite{BXZ1} \cite{Me1} \cite{GR1}, the main result of
this section is to present another way to view this basic gerbe.

Let us recall two lemmas about equivariant $K$-theory and equivariant
cohomology of $G$ \cite{Bj2} \cite{BZ1}.

\begin{lemma}
For a compact simple simply-connected Lie group $G$, 
\begin{equation*}
H^*_G(G)\cong (H^*(BT)\otimes H^*(T))^W
\end{equation*}
\end{lemma}

\begin{lemma}
For a compact simple simply-connected Lie group $G$, 
\begin{equation*}
K^*_G(G)\cong (R(T)\otimes K^*(T))^W
\end{equation*}
\end{lemma}

\begin{prop}
For a compact simple simply-connected Lie group $G$, the basic equivariant
gerbe can be viewed as an element of $K^1_G(G)$.
\end{prop}

{\noindent\textit{Proof.\ \ }} By above lemmas, 
\begin{equation*}
H^3_G(G)\cong (H^0(BT)\otimes H^3(T)\oplus H^2(BT)\otimes H^1(T) )^W
\end{equation*}
\begin{equation*}
\subset (R(T)\otimes K^1(T))^W\cong K^1_G(G),
\end{equation*}
here, $H^0(BT)\cong \mathbb{Z}$, $H^2(BT)\cong X^*(T)$ can be viewed as
subset of $R(T)$. $\square$

\section{$K$-Twisted $K$-theory for $SU(N)$}

In this section, we will use our definition of $K$-twisted equivariant $K$
-theory and Grothendieck differentials to do the calculation for $SU(N)$. 

Let us first recall some background of Grothendieck differentials. Let $
A\subset B$ be commutative rings. The algebra of Grothendieck differentials $
\Omega _{B/A}^{\ast }$ \cite{GD} is the differential graded $A$-algebra
constructed as follows:

Let $F$ be the free $B$-module generated by all elements in B, to be clear,
we use $db$ to denote the generator corresponding to $b \in B$, so 
\begin{equation*}
F=\bigoplus _{b\in B} Bdb.
\end{equation*}
and let $I \subset F$ be the $B$-submodule generated by 
\begin{equation*}
\left\{ 
\begin{array}{l}
da,\ \forall a \in A \\ 
d(b_1+b_2)-db_1-db_2, \ \forall b_1, b_2 \in B \\ 
d(b_1b_2)-b_1db_2-b_2db_1,\ \forall b_1, b_2 \in B
\end{array}
\right\},
\end{equation*}
we then get the quotient $B$-module 
\begin{equation*}
\Omega _{B/A}=F/I.
\end{equation*}

Let $\Omega ^0_{B/A}=B$, $\Omega ^1_{B/A}=\Omega _{B/A}$, and $\Omega
^p_{B/A}= \Lambda ^p_B \Omega _{B/A}$. There is a differential: $d: \Omega
^p_{B/A} \to \Omega ^{p+1}_{B/A}$, which maps $b \in B$ to $db$, then 
\begin{equation*}
\Omega ^*_{B/A}=\bigoplus _{p=0}^{\infty} \Omega ^p_{B/A}
\end{equation*}
is the differential graded algebra of Grothendieck differentials of $B$ over 
$A$. It is the generalization of the algebra of differentials on affine
spaces, for example, if $B=A[x_1, \cdots, x_n]$, then $\Omega ^p_{A[x_1,
\cdots, x_n]/A}=\oplus _{i_1<i_2< \cdots <i_p} A[x_1, \cdots,
x_n]dx_{i_1}\wedge \cdots \wedge dx_{i_p}$.

For any representation $\rho: G\to GL(V)$, it defines a vector bundle over
the suspension of $G$, which is $G$-equivariant, so it defines an element $d\rho$ of $K^1_G(G)$. The main result in \cite{BZ1} is this defines an
isomorphism $\Omega ^*_{R(G)/\mathbb{Z}} \cong K_G^*(G)$, when $\pi _1(G)$ is torsion free.

This result applies to the case of a torus $T$. In terms of Grothedieck
differentials, for any character $\chi _i$ of $T$, $d\chi _i =\chi _i \eta _i$,
where $\eta _i$ is the $K$-theory element or cohomology element of $T$
constructed in the previous section, or in other words, $\eta _i=\frac {
d\chi _i}{\chi _i}$.

In the case $G=SU(N)$, if we let $\rho _i$ be the $i$-th elementary symmetric polynomial in $\chi _1, \chi _2, \cdots, \chi _N$, then $R(G)=\Z [\rho _1, \rho _2, \cdots, \rho _{N-1}]$, and the equivariant $K$-theory is $K^*_G(G)=\wedge_{R(G)}(d\rho _1, d\rho _2, \cdots, d\rho _{N-1})$. 

\begin{prop}
For $SU(N)$, let $\delta $ be the basic gerbe, than as an element of $K^1_{SU(N)}(SU(N)$, is $$\delta=\sum \chi _i \eta _i$$
$$n\delta =\sum \chi _i ^n\eta _i$$
\end{prop}

Let $\alpha =n\delta$, now we are going to calculate the $K$-twisted $K$-theory $^{\alpha}K^*_G(G)$ for $G=SU(N)$, we need a lemma.

\begin {lemma}
Let $\alpha =\sum x_i^ndx_i$ $n\ge 0$, then the following complex is exact except at the last spot:
$$
0\to \Z [x_1, x_2, \cdots, x_N] \overset{\wedge \alpha}{\to}\oplus \Z [x_1, x_2, \cdots, x_N]dx_i \overset{\wedge \alpha}{\to}
$$
$$
\cdots\overset{\wedge \alpha}{\to} \Z [x_1, x_2, \cdots, x_N]dx_1\cdots dx_N\overset{\wedge \alpha}{\to}0
$$
\end{lemma}

Now it is a standard calculation to get $K$-twisted equivariant $K$-theory, in particular,

\begin{Th}
Let $\alpha =(N+k)\delta$, then $^{\alpha }K^N_{SU(N)}(SU(N))$ is the Verlinde algebra $V_k$ of $SU(N)$ at level $k$.
\end{Th}

{\it Proof } By the above lemma and taking $W$-invariants, the non-trivial term of $^{\alpha }K^*_{SU(N)}(SU(N))$ only appears in degree $N$. If $\alpha =a_i d\rho _i$ (These $a_i \in R(SU(N))$ are classical functions, for example see \cite {GV1}), then the $K$-twisted equivariant $K$-theory is $^{\alpha }K^N_{SU(N)}(SU(N))=R(SU(N))d\rho_1d\rho_2 \cdots d\rho_{N-1}/(a_1, a_2, \cdots a_{N-1})d\rho_1d\rho_2 \cdots d\rho_{N-1}\cong R(SU(N))/(a_1, a_2, \cdots a_{N-1})$, which is the Verlinde algebra at level $k$.


\begin{thebibliography}{99}
\bibitem{Fd1} D. Freed, The Verlinde algebra is twisted equivariant
K-theory, math.RT/0101038
\end{thebibliography}
\end{document}